# Solving the Stieltjes Integral Equation in Explicit Form

Peter C. Schuur

Department of Industrial Engineering and Business Information Systems,
University of Twente, Enschede, the Netherlands
p.c.schuur@utwente.nl

## Abstract

Due to its convolution nature, the Stieltjes integral equation can be diagonalized by Mellin transform. Several explicit resolvent kernels were obtained over the years, all of convolution type. The conditions on the given function under which these convolution kernels are able to solve the equation, are rather restrictive. Purpose of this paper is to solve the Stieltjes integral equation - in explicit form - under more general conditions than has been done so far. In fact, we merely bestow upon the given function the same integrability as upon the unknown function. To solve the equation under this mild condition, we construct a new explicit resolvent kernel. For the solutions obtained, we derive interesting growth properties. The new kernel demonstrates that combining known convolution kernels may well lead to a non-convolution kernel that is more effective.

## 1. Introduction

This paper considers the Stieltjes equation

$$f(x) = g(x) + \lambda \int_0^\infty \frac{1}{x+y} f(y)\, dy \qquad x > 0 \qquad (1)$$

In this integral equation, $g: (0, \infty) \to \mathbb{C}$ is given and $f: (0, \infty) \to \mathbb{C}$ is asked. Furthermore, $\lambda \in \mathbb{C}$ is a given parameter. For notational convenience, one writes $\lambda = \sin(\pi\alpha)/\pi$, where $-\frac{1}{2} \leq \operatorname{Re} \alpha \leq \frac{1}{2}$.

The equation has a long history that started with Stieltjes (1894) [9] who employed the associated integral operator in his work on continued fractions. As for more recent applications, the Stieltjes equation is used by Atkinson and Contogouris (1965) [1] in their study on the N/D decomposition of two-particle partial-waves amplitudes.

The homogeneous equation ($g \equiv 0$) was considered by Weyl (1908) [10], Carleman (1923) [4] and Hyslop (1924) [7], until it was finally solved by Hardy and Titchmarsh (1929) [6], albeit under the condition $\lambda > 0$. They use techniques from real analysis, such as the mean value theorem. Widder (1946) [11] (pp 387-389) reconsiders their proof. Although he also focuses on the case $\lambda > 0$, he introduces estimation techniques that can be readily extended to cover the general case $\lambda \in \mathbb{C}$. This leads to the following solution for the homogeneous equation:

   (i) If $\operatorname{Re} \lambda \leq 0$, then $f \equiv 0$ (2.1)

   (ii) If $\operatorname{Re} \lambda > 0$ and $\lambda \neq 1/\pi$, then $f(x) = Ax^{-\alpha} + Bx^{\alpha-1}$ (2.2)

   (iii) If $\lambda = 1/\pi$, then $f(x) = x^{-1/2}(A + B \ln(x))$. (2.3)

Here, $A$ and $B$ are constants.

The kernel of the integral equation (1) is a homogeneous function of the two variables of degree minus one. Equations of this type have been extensively studied by Dixon (1924) [5] who obtains solutions using the Mellin transform. For the inhomogeneous Stieltjes equation (1), Dixon finds three solutions ([5] p 214), each of the form

$$f(x) = g(x) + \lambda \int_0^\infty R_i(x, y; \alpha) g(y)\, dy \qquad x > 0 \qquad (3)$$

where the three resolvent kernels $R_i(x, y; \alpha), i = 1, 2, 3$ are explicitly given. For $Re\,\lambda \leq 0$ only the kernel $R_1$ is effective. For $Re\,\lambda > 0$, each kernel solves equation (1) for a specific set of functions $g$ and for a specific value range of $\lambda$.

To the best of our knowledge, at present, the state of the art for the case $Re\,\lambda > 0$ is as follows. There exists a collection of resolvent kernels such that for different functions $g$ different kernels may have to be used to ensure convergence. In addition, the total set of functions $g$ for which equation (1) can be solved, is rather limited. To illustrate, the kernels found so far will not work for $g(x) = 1/(\ln(2 + x))^2$.

Our paper intends to mend these shortcomings. We present a new family of explicit resolvent kernels for the case $Re\,\lambda > 0$. Each member of this family is able to solve equation (1) under the mild condition that $g(x)$ should be integrable over $(0, \infty)$ with the weight $1/(1 + x)$. Intuitively, this choice seems natural, since the unknown function $f$ also has this property in order to ensure convergence in (1). We prove that, remarkably, the kernel $R_1$ found by Dixon still works for the latter set of functions $g$ in case $Re\,\lambda < 0$.

The organization of the paper is as follows. Section 2 highlights relevant properties of the three resolvent kernels found by Dixon. In Section 3 we solve the Stieltjes equation for quite general $g(x)$ with the help of a new, explicit resolvent kernel. Section 4 establishes growth properties of the solutions. After a concise discussion in Section 5, our conclusion is given in Section 6.

## 2. Some properties of the resolvent kernels found by Dixon

Key to solving the inhomogeneous Stieltjes equation (1) is to write it as

$$f = g + \lambda f * h \tag{4}$$

where $h(x) = 1/(1 + x)$. Here, the convolution product of two functions $f, g: (0, \infty) \to \mathbb{C}$ is defined as

$$(f * g)(x) = \int_0^\infty f(y) g\left(\frac{x}{y}\right) \frac{1}{y} dy$$

for every $x$ for which the integral exists. Closely related is the Mellin transform

$$(Mf)(s) = \int_0^\infty x^{s-1} f(x)\, dx \quad \operatorname{Re} s = k$$

Specifically, if $f(x)$ and $g(x)$ are integrable over $(0, \infty)$ with the weight $x^{k-1}$, then so is $f * g$ and $(M(f * g))(s) = (Mf)(s)(Mg)(s)$ $\quad \operatorname{Re} s = k$

Diagonalizing (4) by Mellin transform one finds for $0 < k < 1$

$$(Mf)(s) = (Mg)(s) + \lambda (Mf)(s)\, \pi / \sin(\pi s)$$

Solving for $(Mf)(s)$ and applying the inverse Mellin transform, Dixon [5] obtains the three explicit resolvent kernels $R_i(x, y; \alpha), i = 1, 2, 3$ mentioned in the introduction (see (3)). In the next three sections we analyze these.

## 2.1 Auxiliary functions

The three resolvent kernels found in [5] are *convolution* kernels, i.e., (3) can be rewritten as $f = g + \lambda r_i(.\,; \alpha) * g$ where the functions $r_i(x; \alpha)$ are given below.

For $x > 0, x \neq 1$:

$$r_1(x; \alpha) = \begin{cases} \dfrac{1}{\cos(\pi \alpha)} \dfrac{x^{-\alpha} - x^{\alpha+1}}{1 - x^2}, & -\tfrac{1}{2} \leq \operatorname{Re} \alpha < \tfrac{1}{2} \quad \alpha \neq -\tfrac{1}{2} \\[6pt] \dfrac{-2}{\pi} \dfrac{x^{1/2}}{1 - x^2} \ln(x), & \alpha = -\tfrac{1}{2} \end{cases}$$

$$r_2(x; \alpha) = -\frac{1}{x} r_1(x; -\alpha) \quad 0 < \operatorname{Re} \alpha \leq \tfrac{1}{2}$$

$$r_3(x; \alpha) = x^2 r_2(x; \alpha) \quad 0 < \operatorname{Re} \alpha \leq \tfrac{1}{2}$$

For the above values of $\alpha$, contour integration yields the fundamental property:

$$r_i(x; \alpha) = \frac{1}{1 + x} + \frac{\sin(\pi \alpha)}{\pi} \int_0^\infty \frac{1}{x + y} r_i(y; \alpha)\, dy \tag{5}$$

So, $r_i(x; \alpha)$ satisfies the Stieltjes equation for the special case that $g(x) = 1/(1 + x)$. For future reference, note that $\int_0^\infty \max(1, x^{-1}) |r_1(x; \alpha)|\, dx < \infty \quad -\tfrac{1}{2} \leq \operatorname{Re} \alpha < 0 \tag{6}$

## 2.2 Resolvent kernels known so far

Using the above auxiliary functions, the three resolvent kernels found by Dixon [5] can be written as follows:

For $x, y > 0, x \neq y$ $\quad R_i(x, y; \alpha) \equiv \dfrac{1}{y} r_i\left(\dfrac{x}{y}; \alpha\right) \quad i = 1, 2, 3$ $\hfill (7)$

Using (7) we can rewrite (5) (for the same $\alpha$-values) as:

$$R_i(x, z; \alpha) = \frac{1}{x+z} + \frac{\sin(\pi\alpha)}{\pi} \int_0^\infty \frac{1}{x+y} R_i(y, z; \alpha) dy \qquad (8)$$

Note that $\quad R_1(x, y; \alpha) = R_1(y, x; \alpha) = \frac{1}{y} r_1\left(\frac{x}{y}; \alpha\right) = \frac{1}{x} r_1\left(\frac{y}{x}; \alpha\right)$ (symmetric kernel)

$$R_2(x, y; \alpha) = R_3(y, x; \alpha) = -\frac{y}{x} R_1(x, y; -\alpha) \qquad 0 < Re\ \alpha \leq \tfrac{1}{2} \qquad (9)$$

## 2.3 Regions of convergence

In [5], regions of convergence for the three resolvent kernels in (3) are not discussed in detail. Let us shed light on this. For $0 < k < 1$ we introduce the B-space $E_k = \{f: (0, \infty) \to \mathbb{C} \mid \|f\|_{E_k} \equiv \int_0^\infty x^{k-1} |f(x)| dx < \infty\}$. Then the solutions of (1) - when we confine ourselves to $E_k$ - are characterized as follows.

**Proposition 1:** Let $0 < k < 1$. Let $\lambda = \sin(\pi\alpha)/\pi$, with $|Re\ \alpha| \leq \tfrac{1}{2}, Re\ \alpha \neq k$ and $Re\ \alpha \neq 1 - k$. Let $g \in E_k$ be given. Then $f \in E_k$ satisfies the inhomogeneous Stieltjes equation (1) iff $f$ is given by: $f(x) = g(x) + \lambda \int_0^\infty R(x, y; k, \alpha) g(y)\, dy \quad x > 0,$ where

$$R(x, y; k, \alpha) \equiv \begin{cases} R_1(x, y; \alpha) \text{ for } -\tfrac{1}{2} \leq Re\ \alpha \leq \max(0, Re\ \alpha) < k < \min(1, 1 - Re\ \alpha) \\ R_2(x, y; \alpha) \text{ for } \tfrac{1}{2} \leq 1 - Re\ \alpha < k < 1 \\ R_3(x, y; \alpha) \text{ for } 0 < k < Re\ \alpha \leq \tfrac{1}{2} \end{cases}$$

**Proof:** Since the nontrivial solutions for the homogeneous equation - given by (2.2) and (2.3) - do not belong to $E_k$, we know that any solution $f \in E_k$ found for (1) is unique.

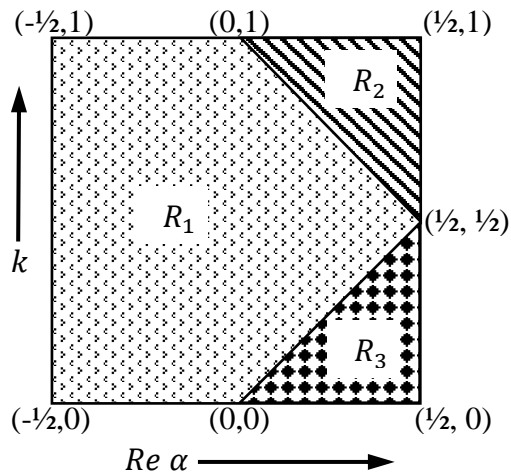

Figure 1: Kernel activity ranges

Note that the above convergence regions for the resolvent kernels $R_i(x, y; \alpha)$ are chosen such that the corresponding functions $r_i(.; \alpha)$ belong to $E_k$.

It remains to be shown that $f = g + \lambda r_i(.;\alpha) * g$ satisfies the Stieltjes equation (4), i.e., $f = g + \lambda f * h$ with $h(x) = 1/(1+x)$. From (5) we know: $r_i(.;\alpha) = h + \lambda r_i(.;\alpha) * h$. Hence, by Fubini: $r_i(.;\alpha) * g = h * g + \lambda(r_i(.;\alpha) * h) * g = (g + \lambda r_i(.;\alpha) * g) * h$.

Consequently, $g + \lambda r_i(.;\alpha) * g = g + \lambda(g + \lambda r_i(.;\alpha) * g) * h$.

In other words, $f = g + \lambda r_i(.;\alpha) * g$ satisfies the Stieltjes equation (4). ∎

The above convergence regions for the kernels $R_i(x,y;\alpha)$ are visualized in Fig.1.

## 3. Constructing new resolvent kernels for $0 < Re\ \alpha \leq \frac{1}{2}$

Using the three resolvent kernels known so far in order to solve the Stieltjes equation (1) may be rather cumbersome. Suppose one wants to solve (1) for a given $g$ and $\lambda$ by applying Proposition 1. Then, one first has to detect a suitable value of $k$ such that $g \in E_k$. Should $Re\ \lambda$ be positive, then one needs the corresponding coordinates $(Re\ \alpha, k)$ to identify which of the three kernels is effective (see Fig. 1).

The following questions come to mind: (A) Is there some overarching space $E$ containing all the $E_k$ such that (1) is solvable for $g \in E$ and $\lambda \in \mathbb{C}$? (B) If so, then specifically for $Re\ \lambda > 0$, can the solution be given in terms of one single resolvent kernel?

Remarkably, the answers to both questions are affirmative, albeit that for question (A) we must confine ourselves to nonzero $Re\ \lambda$.

When searching for $E$, notice that any solution $f(x)$ of (1) is integrable over $(0,\infty)$ with the weight $1/(1+x)$. Thus, for $g$ to guarantee a solution $f$, it should at least be locally integrable on $(0,\infty)$. Moreover, $g(x)$ and $f(x)$ behave similarly, as $x \to \infty$. Hence, it seems natural to bestow upon $g$ the same integrability as upon $f$. Thus, we propose to solve (1) within the B-space $E = \{f: (0,\infty) \to \mathbb{C} \mid \|f\|_E \equiv \int_0^\infty (1+x)^{-1}|f(x)|\,dx < \infty\}$.

As for question (B), notice that we already have the two resolvent kernels $R_2$ and $R_3$. Both have a limited region of convergence. To extend both regions, we propose for $0 < Re\ \alpha \leq \frac{1}{2}$ a new resolvent kernel $R_{23}$ that uses the kernels $R_2$ and $R_3$ as building blocks:

$$R_{23}(x,y;\alpha) \equiv \varphi_1(y)\ R_2(x,y;\alpha) + \varphi_2(y)\ R_3(x,y;\alpha) \qquad (10)$$

where $\varphi_i: (0,\infty) \to (0,1)$, $i = 1, 2$, are continuous functions such that $\varphi_1(y) + \varphi_2(y) = 1$. Thus, for fixed $y$, the new kernel is a convex combination of $R_2$ and $R_3$. It inherits helpful properties like (8). On the other hand, by imposing suitable growth conditions on the functions $\varphi_i$, convergence can be improved considerably. Note that, unlike $R_2$ and $R_3$, the new kernel is not a convolution kernel. Surprisingly, in case $Re\ \lambda < 0$, the convolution kernel $R_1$ found by Dixon [5] still works within $E$. The following theorem specifies our findings.

**Theorem 1:** Let $\lambda \in \mathbb{C}$ and $Re\ \lambda \neq 0$. Write $\lambda = \sin(\pi\alpha)/\pi$, $0 < |Re\ \alpha| \leq \frac{1}{2}$. Let $g \in E$ be given, i.e., $g(x)$ is integrable over $(0, \infty)$ with the weight $1/(1+x)$. Then $f \in E$ satisfies the inhomogeneous Stieltjes equation $f(x) = g(x) + \lambda \int_0^\infty (x+y)^{-1} f(y)\, dy$, $x > 0$, iff $f$ has the following form:

(i) If $Re\ \lambda < 0$ and $\lambda \neq -1/\pi$

$$f(x) = g(x) + \frac{\tan(\pi\alpha)}{\pi} x \int_0^\infty \frac{(y/x)^{-\alpha} - (y/x)^{\alpha+1}}{x^2 - y^2} g(y)\, dy$$

(ii) If $\lambda = -1/\pi$

$$f(x) = g(x) + \frac{2}{\pi^2} x \int_0^\infty \frac{1}{x^2 - y^2} (y/x)^{1/2} \ln(y/x) g(y)\, dy$$

(iii) If $Re\ \lambda > 0$ and $\lambda \neq 1/\pi$

$$f(x) = g(x) - \frac{\tan(\pi\alpha)}{\pi} x \int_0^\infty \left(\varphi_1(y)(y/x) + \varphi_2(y)(x/y)\right) \frac{(y/x)^\alpha - (y/x)^{1-\alpha}}{x^2 - y^2} g(y)\, dy$$
$$+ Ax^{-\alpha} + Bx^{\alpha-1},$$

where $A$ and $B$ are constants.
Here $\varphi_1(y)$ and $\varphi_2(y)$ are continuous positive functions of $y > 0$ such that:
(a) $\varphi_1(y) + \varphi_2(y) = 1$, (b) $\varphi_1(y) = O\left(\frac{1}{y}\right)$ $y \to \infty$, and (c) $\varphi_1(y) = 1 + O(y)$ $y \downarrow 0$

(iv) If $\lambda = 1/\pi$

$$f(x) = g(x) + \frac{2}{\pi^2} x \int_0^\infty \left(\varphi_1(y)(y/x) + \varphi_2(y)(x/y)\right) \frac{1}{x^2 - y^2} (y/x)^{1/2} \ln(y/x) g(y)\, dy$$
$$+ x^{-1/2}(A + B \ln(x)),$$

with $\varphi_1(y)$ and $\varphi_2(y)$ as above.

**Proof:** Let $R(x, y; \alpha) = \begin{cases} R_1(x, y; \alpha) & -½ \leq Re\ \alpha < 0 \\ R_{23}(x, y; \alpha) & 0 < Re\ \alpha \leq ½ \end{cases}$

Since the solutions (2) of the homogeneous equation all belong to $E$, we are done if we can prove that $f(x) = g(x) + \lambda \int_0^\infty R(x,y;\alpha) g(y) dy$ is a particular solution of (1) for $\operatorname{Re} \lambda \neq 0$

Let $0 < |\operatorname{Re} \alpha| \leq \frac{1}{2}$. Let $g \in E$. We show the following:

((1)) $\int_0^\infty |R(x,y;\alpha)| |g(y)| dy < \infty \quad x > 0$

((2)) Let $(Rg)(x) \equiv \int_0^\infty R(x,y;\alpha) g(y) dy$, then $Rg \in E$

((3)) $f = g + \lambda Rg$ satisfies (1) for $\operatorname{Re}\lambda \neq 0$

**Ad ((1))** Deploying (9) we find that $R_{23}(x,y;\alpha) = -\Phi(x,y) R_1(x,y;-\alpha)$
for $0 < \operatorname{Re} \alpha \leq \frac{1}{2}$ with $\Phi(x,y) \equiv \varphi_1(y)(y/x) + \varphi_2(y)(x/y)$.

The properties imposed on $\varphi_i$ guarantee that for fixed $x > 0$, the function $\Phi(x,y)$ is bounded on $y > 0$. Hence, it is sufficient to prove that for $0 < \operatorname{Re} \alpha \leq \frac{1}{2}$

$\int_0^\infty |R_1(x,y;-\alpha)| |g(y)| dy < \infty \quad x > 0$.

The latter follows from $\int_0^\infty |R_1(x,y;-\alpha)| |g(y)| dy \leq M_{\alpha,x} \|g\|_E$

where $M_{\alpha,x} = \sup_{y>0} (1+y) y^{-1} |r_1(x/y;-\alpha)| < \infty$.

**Ad ((2))** We prove more, namely $\int_0^\infty \frac{1}{1+x} \int_0^\infty |R(x,y;\alpha)| |g(y)| dy dx < \infty$ \hfill (11)

Interchanging the integrals, we are faced with proving that $\int_0^\infty \frac{|g(y)|}{1+y} K(y;\alpha) dy < \infty$,

where $K(y;\alpha) \equiv \int_0^\infty \frac{1+y}{1+x} |R(x,y;\alpha)| dx$. Below, we show that $S_\alpha \equiv \sup_{y>0} K(y;\alpha) < \infty$, so that $\int_0^\infty \frac{|g(y)|}{1+y} K(y;\alpha) dy \leq S_\alpha \|g\|_E < \infty$.

Let $-\frac{1}{2} \leq \operatorname{Re} \alpha < 0$.

Then $K(y;\alpha) = \int_0^\infty \frac{1+y}{1+x} \frac{1}{x} |r_1(\frac{y}{x};\alpha)| dx = \int_0^\infty \frac{1+y}{u+y} |r_1(u;\alpha)| du \leq$

$\int_0^\infty \max(1, u^{-1}) |r_1(u;\alpha)| du < \infty$ in view of (6), so that $S_\alpha < \infty$.

Next, let $0 < \operatorname{Re} \alpha \leq \frac{1}{2}$.

Then $K(y;\alpha) = \int_0^\infty \frac{1+y}{1+x} \Phi(x,y) \frac{1}{y} |r_1(\frac{x}{y};-\alpha)| dx = \int_0^\infty L(u,y) |r_1(u;-\alpha)| du$ with

$L(u,y) \equiv \frac{1+y}{1+uy} \Phi(uy,y) = \frac{1+y}{1+uy} (\varphi_1(y)(\frac{1}{u}) + \varphi_2(y) u)$.

Looking at $L(u, y)$, for a fixed value of $y$, we notice that $L(u, y) \approx \frac{1}{u}(1 + y)\varphi_1(y)$ $u \downarrow 0$

and $L(u, y) \approx \frac{1}{y}(1 + y)\varphi_2(y)$ $u \to \infty$. This suggests that $L(u, y)$ may be $O\left(\frac{1}{u}\right)$ for small $u$ and $O(1)$ for large $u$. This motivates us to examine the following two cases:

**Case A** Let $u \leq y\,\varphi_1(y)/\varphi_2(y)$. Then

$uL(u, y) = \frac{1+y}{1+uy}(\varphi_1(y) + \varphi_2(y)u^2) \leq \frac{1+y}{1+uy}\varphi_1(y)(1 + yu) = (1 + y)\varphi_1(y) \leq M_1$ where $M_1 \equiv \sup_{y>0}(1 + y)\varphi_1(y) < \infty$

**Case B** Let $u > y\,\varphi_1(y)/\varphi_2(y)$. Then

$L(u, y) = \frac{1+y}{1+uy}\left(\varphi_1(y)\left(\frac{1}{u}\right) + \varphi_2(y)u\right) < \frac{1+y}{1+uy}\varphi_2(y)\left(\frac{1}{y}\right)(1 + yu) = \frac{1+y}{y}\varphi_2(y) \leq M_2$

where $M_2 \equiv \sup_{y>0}(1 + y)y^{-1}\varphi_2(y) < \infty$

Together, the two cases yield: $L(u, y) \leq M_3 \max(1, u^{-1})$ with $M_3 = \max(M_1, M_2)$.

Hence $K(y; \alpha) = \int_0^\infty L(u, y)\,|r_1(u; -\alpha)|du \leq \int_0^\infty M_3 \max(1, u^{-1})\,|r_1(u; -\alpha)|du < \infty$ in view of (6), so that $S_\alpha < \infty$.

**Ad ((3))** We first show that for $Re\ \lambda \neq 0$:

$$R(x, z; \alpha) = \frac{1}{x + z} + \lambda \int_0^\infty \frac{1}{x + y} R(y, z; \alpha)dy \qquad (12)$$

For $Re\ \lambda < 0$, this fundamental property is already known from (8).

Let us consider $Re\ \lambda > 0$. Then, using (8) for $R_2$ and $R_3$, we obtain $R(x, z; \alpha) =$

$\varphi_1(z)\left(\frac{1}{x+z} + \lambda \int_0^\infty \frac{1}{x+y} R_2(y, z; \alpha)dy\right) + \varphi_2(z)\left(\frac{1}{x+z} + \lambda \int_0^\infty \frac{1}{x+y} R_3(y, z; \alpha)dy\right) = \frac{1}{x+z} +$

$\lambda \int_0^\infty \frac{1}{x+y}(\varphi_1(z)\ R_2(y, z; \alpha) + \varphi_2(z)\ R_3(y, z; \alpha))dy = \frac{1}{x+z} + \lambda \int_0^\infty \frac{1}{x+y} R(y, z; \alpha)dy$ ∎

From (12) we see

$(Rg)(x) = \int_0^\infty R(x, z; \alpha)g(z)\,dz = \int_0^\infty \frac{1}{x+z}g(z)\,dz + \lambda \int_0^\infty \int_0^\infty \frac{1}{x+y} R(y, z; \alpha)g(z)\,dydz$

Because of (11), we have that $\int_0^\infty \frac{1}{x+y} \int_0^\infty |R(y, z; \alpha)|\,|g(z)|dzdy < \infty$. Hence. by Fubini,

we may interchange the integrals, yielding: $(Rg)(x) = \int_0^\infty \frac{1}{x+y}g(y)\,dy +$

$+ \lambda \int_0^\infty \frac{1}{x+y} \int_0^\infty R(y, z; \alpha)\,g(z)\,dzdy = \int_0^\infty \frac{1}{x+y}(g + \lambda Rg)(y)\,dy$

Hence $(g + \lambda Rg)(x) = g(x) + \lambda \int_0^\infty \frac{1}{x+y}(g + \lambda Rg)(y)\,dy$.

Thus $f = g + \lambda Rg$ satisfies the Stieltjes equation $f(x) = g(x) + \lambda \int_0^\infty \frac{1}{x+y} f(y)\, dy$ QED

Note that the class $\Phi$ of functions $\varphi_1$ satisfying the above requirements is far from empty. It contains interesting functions such as $\varphi_1(y) = 1/(1 + y^m)\ m \geq 1$. In addition, if $\varphi_1^A, \varphi_1^B \in \Phi$, then the same holds for their product and for any convex combination. In practice, for a given $g$ and $Re\ \lambda > 0$, one may choose a specific member $R^+(x, y; \alpha)$ of the family of resolvent kernels (10) that allows for an efficient numerical integration (cf. Barseghyan (2017) [3]).

The case $Re\ \lambda = 0$ was left untouched by the above theorem. From Proposition 1 we know that for purely imaginary $\lambda$-values the Stieltjes equation (1) is solvable within a dense subspace of $E$ namely any of the $E_k$. Extending this solvability to all of $E$ is not an option as shown below.

**Proposition 2:** Let $\lambda \in \mathbb{C} \setminus \{0\}$ and $Re\ \lambda = 0$. Then there exist $g \in E$ such that the Stieltjes equation $f(x) = g(x) + \lambda \int_0^\infty (x + y)^{-1} f(y)\, dy$ has no solution $f \in E$.

**Proof:** Let us define the Stieltjes operator $S$ by: $D_S = \{f \in E \mid Sf \in E\}$ and $(Sf)(x) = \int_0^\infty (x + y)^{-1} f(y)\, dy$. It is readily verified that $S$ is a densely defined, closed operator (cf. Kato (1980) [8] p 164). From the classical inversion formula by Stieltjes [9] (see also [11] p 340), we infer that the operator $S$ is injective. Hence, the only eigenfunctions of $S$ are given by (2.2) and (2.3). Consequently, the point spectrum of $S$ is given by $\sigma_p(S) = \{\pi/\sin(\pi\alpha) \mid 0 < Re\ \alpha \leq \frac{1}{2}\} = \{\mu \in \mathbb{C} \mid Re\ \mu > 0\}$.

Now, let $\lambda \in \mathbb{C} \setminus \{0\}$ and $Re\ \lambda = 0$. Suppose that (1) has a solution $f \in E$ for all $g \in E$. Then $I - \lambda S$ is surjective. Put $\mu = 1/\lambda$. Since $\mu \notin \sigma_p(S)$ we have that $\mu - S$ is injective. Hence, $\mu - S = (1/\lambda)(I - \lambda S)$ is bijective, so that $\mu \in \rho(S)$, the resolvent set of $S$. Since the resolvent set of a closed operator is an open set in the complex plane (see [8] p 174), we can find an environment of $\mu$ that is entirely within $\rho(S)$. Bearing in mind that $Re\ \mu = 0$ this would mean that there is a $\mu^1 \in \rho(S)$ with $Re\ \mu^1 > 0$. Thus $\rho(S)$ and $\sigma_p(S)$ would have a nonempty intersection, which is absurd.

## 4. Growth properties of the solutions

Theorem 1 solves the Stieltjes equation (1) under mild restrictions on the given function $g$. In practice, more information may be available about $g$. One may wonder to what extent this information is transferred to any solution $f$. As for transferring integrability properties, we saw in Proposition 1 that for each $g \in E_k$, there is a unique solution $f \in E_k$ for almost all $\lambda \in \mathbb{C}$. Another interesting property of $g$ that may be transferrable to $f$ is its asymptotic behavior. In this section we take $g$ to be continuous and impose growth conditions near 0 and $\infty$. We examine to what extent these growth conditions are shared by any solution $f$.

The proof of Proposition 2 introduces the Stieltjes operator $S$ given by: $D_S = \{f \in E \mid Sf \in E\}$ and $(Sf)(x) = \int_0^\infty (x+y)^{-1} f(y)\, dy$. Evidently, if $g$ belongs to $F$ - where $F$ is some subspace of $E$ - then also $f$ belongs to $F$ in case both the Stieltjes operator and the corresponding resolvent operator leave $F$ invariant.

### 4.1 A new subspace invariant under the Stieltjes operator

Let $0 < \varepsilon, \eta < 1$. We impose growth conditions on the known function $g$ by choosing it in the B-space $B_{\varepsilon,\eta}$ of all continuous functions $f: (0, \infty) \to \mathbb{C}$ such that $\|f\|_{\varepsilon,\eta} \equiv \sup_{1 < x < \infty} |x^\varepsilon f(x)| + \sup_{0 < x < 1} |x^\eta f(x)| < \infty$.

The B-space $E_k$ introduced in Section 2.3 is invariant under $S$, as a consequence of the convolution structure of $S$. Moreover, it is readily seen that $S$ is a bounded linear operator when confined to $E_k$. Surprisingly, as the next proposition shows, similar statements hold for $B_{\varepsilon,\eta}$.

**Proposition 3:** Let $0 < \varepsilon, \eta < 1$. Then $B_{\varepsilon,\eta}$ is contained in $D_S$. Let $S_{\varepsilon,\eta}$ denote the operator $S$ confined to $B_{\varepsilon,\eta}$. Then $S_{\varepsilon,\eta}$ is a bounded linear operator from $B_{\varepsilon,\eta}$ to $B_{\varepsilon,\eta}$.

**Proof:** For $f \in B_{\varepsilon,\eta}$ one has:

$\int_0^\infty (x+y)^{-1} |f(y)|\, dy \leq \|f\|_{\varepsilon,\eta} (\int_1^\infty (x+y)^{-1} y^{-\varepsilon}\, dy + \int_0^1 (x+y)^{-1} y^{-\eta}\, dy)$.

For $x \geq 1$ this yields:

$\int_0^\infty (x+y)^{-1} |f(y)|\, dy \leq \|f\|_{\varepsilon,\eta} (\int_0^\infty (x+y)^{-1} y^{-\varepsilon}\, dy + \int_0^1 x^{-1} y^{-\eta}\, dy) =$

$= \|f\|_{\varepsilon,\eta}((\pi/\sin(\pi\varepsilon))\, x^{-\varepsilon} + (1-\eta)^{-1} x^{-1}) \leq \|f\|_{\varepsilon,\eta}((\pi/\sin(\pi\varepsilon)) + (1-\eta)^{-1}) x^{-\varepsilon}$.

For $0 < x < 1$ it yields:

$$\int_0^\infty (x+y)^{-1}|f(y)|\,dy \le \|f\|_{\varepsilon,\eta}\left(\int_1^\infty y^{-1}y^{-\varepsilon}\,dy + \int_0^\infty (x+y)^{-1}y^{-\eta}\,dy\right) =$$

$$= \|f\|_{\varepsilon,\eta}(\varepsilon^{-1} + (\pi/\sin(\pi\eta))x^{-\eta}) \le \|f\|_{\varepsilon,\eta}(\varepsilon^{-1} + (\pi/\sin(\pi\eta)))x^{-\eta}.$$

By dominated convergence, the function $\int_0^\infty (x+y)^{-1}f(y)\,dy$ is continuous on $x > 0$. Thus, $B_{\varepsilon,\eta}$ is a subspace of $D_S$, which is left invariant by the bounded linear operator $S_{\varepsilon,\eta}$.

Invariance of $B_{\varepsilon,\eta}$ under the various resolvent operators is a more delicate matter that we discuss in the next section. Crucial there is the following result due to Atkinson (1979) [2].

**Proposition 4:** Let $\beta \in \mathbb{C}$ such that $-1 < \operatorname{Re} \beta < 2$. Let $\max(0, \operatorname{Re} \beta - 1) < \varepsilon, \eta < \min(1, \operatorname{Re}\beta + 1)$. For $f \in B_{\varepsilon,\eta}$ put

$$(T^\beta f)(x) \equiv x \int_0^\infty (y^2 - x^2)^{-1}\bigl((y/x)^\beta - 1\bigr)f(y)\,dy \quad x > 0$$

Then $T^\beta$ is a bounded linear operator from $B_{\varepsilon,\eta}$ to $B_{\varepsilon,\eta}$.

**Proof:** Putting $s(u) = (u^2 - 1)^{-1}(u^\beta - 1)$, the RHS reads $\int_0^\infty x^{-1}s(y/x)f(y)\,dy$. Note that the restriction on $\varepsilon, \eta$, and $\beta$ implies that $\int_0^\infty \max(u^{-\varepsilon}, u^{-\eta})|s(u)|\,du < \infty$. Clearly

$$\int_0^\infty x^{-1}|s(y/x)f(y)|\,dy \le \|f\|_{\varepsilon,\eta}\left(\int_1^\infty x^{-1}|s(y/x)|y^{-\varepsilon}\,dy + \int_0^1 x^{-1}|s(y/x)|y^{-\eta}\,dy\right).$$

For $x \ge 1$ this yields: $\int_0^\infty x^{-1}|s(y/x)f(y)|\,dy \le$

$$\le \|f\|_{\varepsilon,\eta}\left(\int_0^\infty x^{-1}|s(y/x)|y^{-\varepsilon}\,dy + \int_0^x x^{-1}|s(y/x)|y^{-\max(\varepsilon,\eta)}\,dy\right) =$$

$$= \|f\|_{\varepsilon,\eta}\left(x^{-\varepsilon}\int_0^\infty |s(u)|u^{-\varepsilon}\,du + x^{-\max(\varepsilon,\eta)}\int_0^1 |s(u)|u^{-\max(\varepsilon,\eta)}\,du\right) \le$$

$$\le \|f\|_{\varepsilon,\eta}\left(\int_0^\infty |s(u)|u^{-\varepsilon}\,du + \int_0^1 |s(u)|u^{-\max(\varepsilon,\eta)}\,du\right)x^{-\varepsilon}, \text{ where both integrals are finite.}$$

For $0 < x < 1$ it yields: $\int_0^\infty x^{-1}|s(y/x)f(y)|\,dy \le$

$$\le \|f\|_{\varepsilon,\eta}\left(\int_x^\infty x^{-1}|s(y/x)|y^{-\min(\varepsilon,\eta)}\,dy + \int_0^\infty x^{-1}|s(y/x)|y^{-\eta}\,dy\right) =$$

$$= \|f\|_{\varepsilon,\eta}\left(x^{-\min(\varepsilon,\eta)}\int_1^\infty |s(u)|u^{-\min(\varepsilon,\eta)}\,du + x^{-\eta}\int_0^\infty |s(u)|u^{-\eta}\,du\right) \le$$

$$\le \|f\|_{\varepsilon,\eta}\left(\int_1^\infty |s(u)|u^{-\min(\varepsilon,\eta)}\,du + \int_0^\infty |s(u)|u^{-\eta}\,du\right)x^{-\eta}, \text{ which completes the proof.}$$

## 4.2 Invariance analysis of the various resolvent operators

From Theorem 1 we recall that $f(x) = g(x) + \lambda \int_0^\infty R(x,y;\alpha) g(y) dy$ is a particular solution of (1) for $Re\ \lambda \neq 0$, where $R(x,y;\alpha) = \begin{cases} R_1(x,y;\alpha) & -½ \leq Re\ \alpha < 0 \\ R_{23}(x,y;\alpha) & 0 < Re\ \alpha \leq ½ \end{cases}$ with

$R_{23}(x,y;\alpha) = \varphi_1(y)\ R_2(x,y;\alpha) + \varphi_2(y)\ R_3(x,y;\alpha)$.

Let us examine for which values of $\varepsilon, \eta$ and $\alpha$ the space $B_{\varepsilon,\eta}$ is invariant under $R(x,y;\alpha)$. For $f \in B_{\varepsilon,\eta}$ we write $(R_i^\alpha f)(x) = \int_0^\infty R_i(x,y;\alpha) f(y) dy\ \ i = 1, 2, 3$, and $23$, whenever the integral exists.

The following two propositions show that the $R_i^\alpha$ leave $B_{\varepsilon,\eta}$ invariant for specific parameter values.

**Proposition 5:** The space $B_{\varepsilon,\eta}$ is invariant under

(i) $R_1^\alpha$ for $-½ \leq Re\ \alpha < 0$ and $0 < \varepsilon, \eta < 1$

(ii) $R_2^\alpha$ for $0 < Re\ \alpha \leq ½$ and $1 - Re\ \alpha < \varepsilon, \eta < 1$

(iii) $R_3^\alpha$ for $0 < Re\ \alpha \leq ½$ and $0 < \varepsilon, \eta < Re\ \alpha$

**Proof:** Due to limit transition the $\alpha$-regions may be confined to $-½ \leq Re\ \alpha < 0, \alpha \neq -½$ for $R_1^\alpha$ and to $0 < Re\ \alpha \leq ½, \alpha \neq ½$ for $R_2^\alpha$ and $R_3^\alpha$. Then Proposition 4 yields for $f \in B_{\varepsilon,\eta}$

$R_1^\alpha f = (\cos(\pi\alpha))^{-1}(T^{\alpha+1} - T^{-\alpha}) f \in B_{\varepsilon,\eta}$ for $0 < \varepsilon, \eta < 1$

$R_2^\alpha f = (\cos(\pi\alpha))^{-1}(T^{\alpha+1} - T^{2-\alpha}) f \in B_{\varepsilon,\eta}$ for $1 - Re\ \alpha < \varepsilon, \eta < 1$

$R_3^\alpha f = (\cos(\pi\alpha))^{-1}(T^{\alpha-1} - T^{-\alpha}) f \in B_{\varepsilon,\eta}$ for $0 < \varepsilon, \eta < Re\ \alpha$

**Proposition 6:** Let $0 < Re\ \alpha \leq ½$. Let $\varphi_1$ and $\varphi_2$ be as in Theorem 1. Then the space $B_{\varepsilon,\eta}$ is invariant under $R_{23}^\alpha$ provided that $0 < \varepsilon < Re\ \alpha$ and $1 - Re\ \alpha < \eta < 1$.

**Proof:** Let $f \in B_{\varepsilon,\eta}$ with $0 < \varepsilon < Re\ \alpha$ and $1 - Re\ \alpha < \eta < 1$. Put $f_1(y) = \varphi_1(y) f(y)$ and $f_2(y) = \varphi_2(y) f(y)$. Then $f_1 \in B_{\hat\varepsilon,\eta}$ for $0 < \hat\varepsilon < 1$ since $\varphi_1(y) = O\left(\frac{1}{y}\right)$ as $y \to \infty$.

Also, $f_2 \in B_{\varepsilon,\hat\eta}$ for $0 < \hat\eta < 1$ since $\varphi_2(y) = O(y)$ as $y \downarrow 0$.

In particular, one has $f_1 \in B_{1-½Re\ \alpha,\eta}$ and $f_2 \in B_{\varepsilon,½Re\ \alpha}$.

By Proposition 5: $R_{23}^\alpha f = R_2^\alpha f_1 + R_3^\alpha f_2 \in B_{\min(\varepsilon, 1-½Re\ \alpha),\max(½Re\ \alpha,\eta)} = B_{\varepsilon,\eta}$

In the proof of Proposition 6 we used the evident inclusion relations:

$(i)$ $B_{\varepsilon_2,\eta} \subset B_{\varepsilon_1,\eta}$ if $\varepsilon_1 < \varepsilon_2$  and  $(ii)$ $B_{\varepsilon,\eta_1} \subset B_{\varepsilon,\eta_2}$ if $\eta_1 < \eta_2$

These inclusion relations give rise to the following corollary:

**Corollary 1:** Let $0 < Re\ \alpha \leq \frac{1}{2}$ and $0 < \varepsilon, \eta < 1$. Let $f \in B_{\varepsilon,\eta}$. Then $R_{23}^{\alpha} f \in B_{\tilde{\varepsilon},\tilde{\eta}}$ where

(i) $\tilde{\varepsilon}$ only depends on $\varepsilon$ and (ii) $\tilde{\eta}$ only depends on $\eta$, as follows:

If $0 < \varepsilon < Re\ \alpha$ then $\tilde{\varepsilon} = \varepsilon$, else $\tilde{\varepsilon}$ is any value satisfying $0 < \tilde{\varepsilon} < Re\ \alpha$

If $1 - Re\ \alpha < \eta < 1$ then $\tilde{\eta} = \eta$, else $\tilde{\eta}$ is any value satisfying $1 - Re\ \alpha < \tilde{\eta} < 1$

Fig. 2 illustrates the selection of $\tilde{\varepsilon}$ and $\tilde{\eta}$ for any $\varepsilon$ and $\eta$ in the parameter set. Combining Theorem 1, Proposition 5, and Corollary 1 we arrive at:

**Corollary 2:** Consider the Stieltjes equation (1), with $g \in B_{\varepsilon,\eta}$ given and $f \in E$ asked.

If $Re\lambda < 0$, then (1) has a unique solution. This solution belongs to $B_{\varepsilon,\eta}$.

If $Re\lambda > 0$, then (1) has a particular solution that belongs to $B_{\tilde{\varepsilon},\tilde{\eta}}$, where $\tilde{\varepsilon}$ and $\tilde{\eta}$ are selected according to Fig. 2.

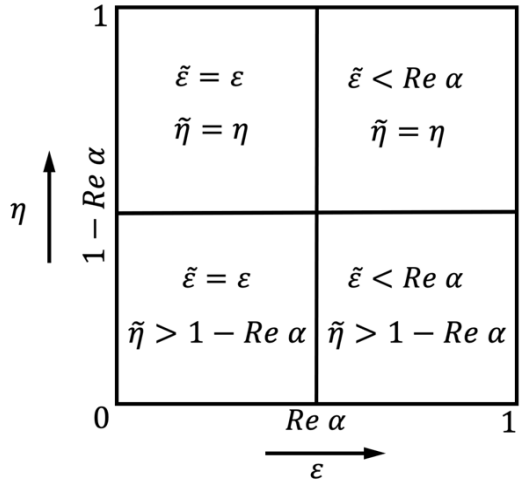

Figure 2: Selection of $\tilde{\varepsilon}$ and $\tilde{\eta}$ for any $\varepsilon$ and $\eta$

**Remark 1:** For $Re\ \lambda = 0$, one can show that (1) with $g \in B_{\varepsilon,\eta}$ is uniquely solvable and that the solution belongs to $B_{\varepsilon,\eta}$

**Remark 2:** Let $Re\ \lambda > 0$ be given. One may wonder: are there concrete instances of $\varepsilon$, $\eta$ and $g \in B_{\varepsilon,\eta}$ such that equation (1) has no solution $f \in B_{\varepsilon,\eta}$? The answer is affirmative. There exists a function $g$ such that: (i) $g \in B_{\varepsilon,\eta}$ for all $0 < \varepsilon, \eta < 1$ and (ii) equation (1) has no solution $f \in B_{\varepsilon,\eta}$ for $\eta < Re\ \alpha < \varepsilon$ and for $\eta < 1 - Re\ \alpha < \varepsilon$.

To see this, take $g(x) = 1/(1 + x)$. Then $g \in B_{\varepsilon,\eta}$ for all $0 < \varepsilon, \eta < 1$. Combining (5) with Theorem 1, we find that the associated solution $f \in E$ of (1) satisfies $f(x) = r_1(x; \alpha) + Ax^{-\alpha} + Bx^{\alpha-1}$ for $x > 0, x \neq 1$. Estimating the RHS, one readily verifies that $f \notin B_{\varepsilon,\eta}$ for $\eta < Re\ \alpha < \varepsilon$ and for $\eta < 1 - Re\ \alpha < \varepsilon$.

## 5. Discussion

The results of this paper lean heavily on three remarkable characteristics of the Stieltjes equation: (i) It can be diagonalized by Mellin transform yielding explicit resolvent kernels. (ii) A convex combination of two specific kernels gives a new explicit resolvent kernel with a large convergence range. (iii) The explicit nature of the solutions unveils intriguing growth properties.

In general, explicit solutions are rare. To illustrate, if we change the lower bound in the Stieltjes equation from 0 to 1, then it can still be diagonalized - by Mehler transform - but the resulting resolvent kernels lack an apparent explicit structure, which makes any form of analysis challenging. Nevertheless, even in the case of an integral equation with a solution that is less accessible, it could be valuable to explore whether a convex combination of resolvent kernels produces a more appropriate kernel.

## 6. Conclusion

This paper establishes that the Stieltjes integral equation (1) can be solved for $Re\ \lambda \neq 0$ under more general conditions on $g$ than has been done so far. We only impose the mild condition that $g \in E$ i.e., $g(x)$ should be integrable over $(0, \infty)$ with the weight $1/(1 + x)$. We prove that, remarkably, the resolvent kernel found by Dixon still works within $E$ in case $Re\ \lambda < 0$. For the case $Re\ \lambda > 0$, the known convolution kernels are not effective on $E$. By combining two of these known convolution kernels, we construct a new family of explicit resolvent kernels. Each member of this family is a non-convolution kernel able to solve equation (1) within $E$ for $Re\ \lambda > 0$. For the solutions obtained in this fashion, we derive interesting growth properties. Our findings suggest that instead of sticking to convolution kernels it may be worthwhile to consider whether combining them into a non-convolution kernel may be more effective.

## Acknowledgements

I am grateful to David Atkinson for introducing the Stieltjes equation to me. Moreover, I thank Erik Thomas and David for stimulating discussions.